\documentclass[10pt,reqno]{amsart}
\usepackage{indentfirst, amssymb,amsmath,amsthm}
\newtheorem*{theoA}{Theorem A}
\newtheorem*{theoB}{Theorem B}
\newtheorem*{theoC}{Theorem C}
\newtheorem*{theoD}{Theorem D}
\newtheorem*{theoE}{Theorem E}
\newtheorem*{theoF}{Theorem F}

\newtheorem{theo}{Theorem}
\newtheorem{lem}{Lemma}
\newtheorem{cor}{Corollary}

\newtheorem{exm}{Example}
\newtheorem{rem}{Remark}

\newcommand{\be}{\begin{equation}}
\newcommand{\ee}{\end{equation}}
\newcommand{\beas}{\begin{eqnarray*}}
\newcommand{\eeas}{\end{eqnarray*}}
\newcommand{\bea}{\begin{eqnarray}}
\newcommand{\eea}{\end{eqnarray}}
\renewcommand{\epsilon}{\varepsilon}
\makeatletter
\renewcommand*\env@matrix[1][*\c@MaxMatrixCols c]{%
\hskip -\arraycolsep
\let\@ifnextchar\new@ifnextchar
\array{#1}}
\makeatother

\numberwithin{equation}{section}
\numberwithin{lem}{section}
\numberwithin{theo}{section}
\numberwithin{cor}{section}
\numberwithin{exm}{section}
\numberwithin{defi}{section}
\numberwithin{rem}{section}

\begin{document}
\title[On transcendental meromorphic solutions...]
{On transcendental meromorphic solutions of certain types of differential equations}
\date{}
\author{Tania Biswas, Sayantan Maity and Abhijit Banerjee}
\address{ Department of Mathematics, University of Kalyani, West Bengal 741235, India.}
\email{taniabiswas2394@gmail.com}
\address{ Department of Mathematics, University of Kalyani, West Bengal 741235, India.}
\email{sayantanmaity100@gmail.com}
\address{ Department of Mathematics, University of Kalyani, West Bengal 741235, India.}
\email{abanerjee\_kal@yahoo.co.in, abanerjeekal@gmail.com}
\renewcommand{\thefootnote}{}
\footnote{2010 {\emph{Mathematics Subject Classification}}: 34A34, 34M05, 30D35. }
\footnote{\emph{Key words and phrases}: Non-linear differential equations, Meromorphic solutions. }
\renewcommand{\thefootnote}{\arabic{footnote}}
\setcounter{footnote}{0}
\begin{abstract} In this paper, for a transcendental meromorphic function $f$ and $a\in \mathbb{C}$, we have exhaustively studied the nature and form of solutions of a new type of non-linear differential equation of the following form which has never been investigated earlier: \beas f^n+af^{n-2}f'+ P_d(z,f) = \sum_{i=1}^{k}p_i(z)e^{\alpha_i(z)},\eeas where $P_d(z,f)$ is differential polynomial of $f$, $p_i$'s and $\alpha_{i}$'s  are non-vanishing
	rational functions and non-constant polynomials respectively. When $a=0$, we have pointed out a major lacuna in a recent result of Xue [Math. Slovaca, 70(1)(2020), 87-94] and rectifying the result, presented the corrected form of the same at a large extent. In addition, our main result is also an improvement of a recent result of Chen-Lian [Bull. Korean Math. Soc., 54(4)(2020), 1061-1073] by rectifying a gap in the proof of the theorem of the same paper. The case $a\neq 0$ has also been manipulated to determine the form of the solutions. We also illustrate a handful number of examples for showing the accuracy of our results.    
	
\end{abstract}
\thanks{Typeset by \AmS -\LaTeX}
\maketitle
\section{Introduction}
Let $\mathbb{C}$ denote the field of complex numbers and $\mathcal{M}(\mathbb{C})$ be the field of meromorphic functions on $\mathbb{C}$. Throughout this paper we consider $f\in \mathcal{M}(\mathbb{C})$.
We assume that the readers are familiar with basic Nevanlinna theory  and usual notations such as proximity function $m(r, f)$, counting function $N(r, f)$, characteristic function $T(r, f)$, first and second main theorem etc (see \cite{hayman_book64}). Recall that, $S(r, f) = o(T(r, f))$ as $r \rightarrow \infty$ outside of a possible exceptional set of finite logarithmic (linear) measure.
A meromorphic function $\alpha$ is called a small function of $f$ if and only if
$T(r, \alpha) = S(r, f)$. The order of a meromorphic function $f$ is denoted by  $\rho(f)$ and defined by $\rho(f) = \displaystyle\limsup_{r \rightarrow \infty} \frac{\log^+\; T(r,f)}{\log\; r}$. For $\alpha \in \mathbb{C} $, the deficiency $\delta (\alpha,f)$ is defined as $ \delta(\alpha,f) = 1- \displaystyle\limsup_{r \rightarrow \infty} \frac{N(r,\alpha;f)}{T(r,f)}.$\par

Nowadays differential equation plays a prominent role in many disciplines and so study on the different features of  differential equation over $\mathbb{C}$ has become an interesting topic. Speculations over the existence of solutions of non-linear differential equation and subsequently finding the exact form of the same are really  challenging problems. In the present paper we wish to contribute in this perspective. To this end, we denote by $P_d(z,f)$ as the non-linear differential polynomial of $f(z)$ of degree $d$ defined by \bea \label{e1.1} P_d(z,f) = \sum_{\lambda \in \Lambda } a_{\lambda} \prod_{i=0}^{n_1}\left(f^{(i)}(z)\right)^{\lambda_i},\eea where $a_{\lambda}$'s are rational function, $\Lambda$ be the index set of non-negative integers with finite cardinality and $\lambda = (\lambda_0, \lambda_1, \lambda_2,\cdots, \lambda_{n_1})$ also $ d:= \deg (P_d(z,f)) = \max\limits_{\lambda \in \Lambda} \left\{\displaystyle \sum_{i=0}^{n_1}\lambda_i\right\}$. \par  From last two decades researchers extensively studied the differential equation of the following form \bea \label{e1.2} f^n + P_d(z,f) = p_1(z) e^{\alpha_1(z)} + p_2(z) e^{\alpha_2(z)},\eea where $P_d(z,f)$ is defined as in (\ref{e1.1}) with some restriction on degree $d$ and $p_1(z)$, $p_2(z)$ are non-zero rational functions and $\alpha_1(z),\alpha_2(z)$ are non-constant polynomials. In this paper also we wish to contribute in this aspect under a more general settings.
\section{ Backgrounds and Main Results} In 2006, on the existence of solution of differential equation, Li-Yang \cite{ping li c c yang_jmaa06} obtained the following result. 
\begin{theoA} \cite{ping li c c yang_jmaa06} Let $n\geq 4$ be an integer. Consider the differential equation (\ref{e1.2}), where $d \leq n-3$, $p_j(z)$ ($j=1,2$) be two non-vanishing polynomials and $\alpha_j(z):=\alpha_jz$ ($j=1,2$) are two non-zero one degree polynomials such that $\frac{\alpha_1}{\alpha_2}$ is not rational. Then the
	equation (\ref{e1.2}) has no transcendental entire solutions. \end{theoA}
In 2011, Li \cite{ping li_jmaa11} removed the extra supposition ``$\frac{\alpha_1}{\alpha_2}$ is not rational" in {\it Theorem A} and established the form of the meromorphic solution. 
\begin{theoB} \cite{ping li_jmaa11} Let $n \geq 2$ be an integer, $P_d(z,f)$ be a differential polynomial in $f(z)$ of degree $d \leq n-2$. Consider the differential equation (\ref{e1.2}), where $d \leq n-2$, $p_j(z)$ ($j=1,2$) be two non-zero constants and $\alpha_j(z):=\alpha_jz$ ($j=1,2$) are two non-zero one degree polynomials such that $\alpha_1 \neq \alpha_2$. If $f(z)$ is a transcendental meromorphic solution of the equation (\ref{e1.2}) and satisfying $N(r, f ) = S(r, f )$, then one of the following holds:\par $(i)$ $f(z) = c_0 + c_1e^{\alpha_1z/n}$;\par $(ii)$ $f(z) = c_0 + c_2e^{\alpha_2z/n}$;\par $(iii)$ $f(z) = c_1e^{\alpha_1z/n} + c_2e^{\alpha_2z/n}$, and $\alpha_1 + \alpha_2 = 0$,\\  where $c_0$ is a small function of $f(z)$ and $c_1, c_2$ are constants satisfying $c^n_1 = p_1$, $c^n_2 =p_2$. \end{theoB}    
In 2013, considering $p_1(z)$, $p_2(z)$ as non-vanishing rational functions and $\alpha_1(z)$, $\alpha_2(z)$ as non-constant polynomials, Liao-Yang-Zhang \cite{liao_ ann fenn13} generalized {\it Theorem B} as follows.  
\begin{theoC} \cite{liao_ ann fenn13} Let $n\geq 3$ be an integer, $P_d(z,f)$ be a differential polynomial in $f(z)$ of degree $d$ with rational functions as its coefficients. Suppose $p_1(z)$, $p_2(z)$ are non-vanishing rational functions and $\alpha_1(z),\alpha_2(z)$ are non-constant polynomials. If $d \leq n-2$ and the differential equation (\ref{e1.2}) admits a meromorphic function solution $f(z)$ with finitely many poles, then $\frac{\alpha'_1}{\alpha'_2}$ is a rational number. Furthermore, only one of the following four cases holds:\par
	$(i)$ $f(z) = q(z)e^{P(z)}$, $\frac{\alpha'_1}{\alpha'_2} = 1$, where $q(z)$ is a rational function and $P(z)$ is a polynomial with $nP' = \alpha'_1 = \alpha'_2$;\par
	$(ii)$ $f(z) = q(z)e^{P(z)}$, either $\frac{\alpha'_1}{\alpha'_2} = \frac{n}{k}$ or $ \frac{k}{n}$, where $q(z)$ is a rational function, $k$ is an integer with $1 \leq k \leq d$ and $P(z)$ is a polynomial with $nP' = \alpha'_1$ or $\alpha'_2$;\par
	$(iii)$ $f$ satisfies the first order linear differential equation $f' = \left(\frac{1}{n}\frac{p'_1}{p_1} + \frac{1}{n} \alpha'_1\right)f + \varphi$ and $\frac{\alpha'_1}{\alpha'_2} = \frac{n}{n-1}$ or $f$ satisfies the first order linear differential equation $f' = \left(\frac{1}{n}\frac{p'_2}{p_2} + \frac{1}{n} \alpha'_2\right)f + \varphi$ and $\frac{\alpha'_1}{\alpha'_2} = \frac{n-1}{n}$, where $\varphi$ is a rational function;\par 
	$(iv)$ $f(z) = c_1(z)e^{\beta(z)} + c_2(z)e^{-\beta(z)}$ and $\frac{\alpha'_1}{\alpha'_2} =-1$, where $c_1(z), c_2(z)$ are
	rational functions and $\beta(z)$ is a polynomial with $n\beta' = \alpha'_1$ or $\alpha'_2$.\end{theoC}
In 2018, Zhang \cite{zhang_acta math sinica18} established the following result in this direction.
\begin{theoD} \cite{zhang_acta math sinica18} Under the same assumption as in {\em Theorem C} if $n\geq 4$ be an integer, $d \leq n-3$, and the complex differential equation (\ref{e1.2}) admits a transcendental meromorphic function solution $f$ with finitely many poles, then $\frac{\alpha'_1}{\alpha'_2}$ is a rational number and $f(z)$ must be of the following form {\em :} \beas f(z) = q(z)e^{P(z)}, \eeas where $q(z)$ is a non-vanishing rational function and $P(z)$ is a non-constant
	polynomial. Moreover, only one of the following two cases holds:\par
	$(i)$ $\frac{\alpha'_1}{\alpha'_2} = 1$, $P_d(z,f) \equiv 0$ and $nP' = \alpha'_1 = \alpha'_2$;\par
	$(ii)$ $\frac{\alpha'_1}{\alpha'_2} = \frac{t}{n}$, where $t$ is an integer satisfying $1 \leq t < d$, $P_d(z,f) \equiv p_1(z)e^{\alpha_1(z)}$ and $nP' = \alpha'_2$ or $\frac{\alpha'_1}{\alpha'_2} = \frac{n}{t}$, where $t$ is an integer satisfying $1 \leq t < d$, $P_d(z,f) \equiv p_2(z)e^{\alpha_2(z)}$ and $nP' = \alpha'_1$.
\end{theoD}

For more results related to this area readers can see \cite{chen_j comput anal appl18, laine_book93, liao_aust math soc14, lu lu_result math19, liao_acta math sinica17, rong xu_mathematics19}. We note that the right hand side of (\ref{e1.2}) consisting of two exponential terms. So investigations on the case when the right hand side of (\ref{e1.2}) contains three or more exponential terms attracts the researchers. In this perspective, in 2020, Xue \cite{xue_math slovaca20} found the following result:

\begin{theoE} \cite{xue_math slovaca20}
	Let $n \geq 2$ and $P_d(z,f)$ be an algebraic differential polynomial in $f(z)$ of degree $d\leq n-1$. Suppose that $p_j$, $\alpha_j$ are non-zero constants for  $j = 1, 2, 3$ and $|\alpha_1|> |\alpha_2| > |\alpha_3|$. If $f(z)$ is a transcendental entire solution of the differential equation \beas f^n + P_d(z,f) = p_1 e^{\alpha_1z} + p_2 e^{\alpha_2z} +  p_3 e^{\alpha_3z},\eeas then $f(z) = a_1e^{\alpha_1z/n}$, where $a_1$ is non-zero constant such that $a^n_1
	= p_1$, and $\alpha_j$ are in one line for $j = 1, 2, 3$. \end{theoE} 
\begin{rem}
	Although {\it Theorem E} is a novel approach in this direction but we find that there is a major drawback in the last theorem. The following example shows that {\it Theorem E} does not hold good. \end{rem}

\begin{exm} The function $f(z) = e^z + 1$ is a transcendental entire solution of \beas f^3 + 4f'f + f'-f = e^{3z} + 7e^{2z} +7 e^z.\eeas But $e^z + 1$ is not of the form of $f(z)$ given in {\it Theorem E}. \end{exm}

Also in 2020, Chen-Lian \cite{chen lian_BKMS20} studied on the differential equation with three exponential terms in the right hand side and established the following:

\begin{theoF} \cite{chen lian_BKMS20} Let $n\geq 5$ be an integer, $P_d(z,f)$ be a differential
	polynomial in $f(z)$ of degree $d$ with rational functions as its coefficients. Suppose $p_j(z)$ ($j=1,2,3$) are non-vanishing rational functions and $\alpha_j(z)$ ($j=1,2,3$) are non-constant polynomials such that $\alpha'_{j}(z)$ ($j=1,2,3$) are distinct each other. If $d \leq n-4$ and the differential equation \beas f^n + P_d(z,f) = \sum_{j=1}^{3}p_j(z)e^{\alpha_j(z)}, \eeas admits a transcendental meromorphic function solution $f$ with finitely many poles, then $\frac{\alpha'_1}{\alpha'_2}$, $\frac{\alpha'_2}{\alpha'_3}$ are rational numbers and $f(z)$ must be of the following form: \beas f(z) = q(z)e^{P(z)}, \eeas where $q(z)$ is a non-vanishing rational function and $P(z)$ is a non-constant polynomial. Moreover, there must exist positive integers $l_0, l_1, l_2$ with \{$l_0, l_1, l_2$\} = \{$1,2,3$\} and distinct integers $k_1, k_2$ with $1 \leq k_1, k_2 \leq d$ such that $\alpha'_{l_0}:\alpha'_{l_1}:\alpha'_{l_2} = n: k_1: k_2$, $nP' = \alpha'_{l_0}$ and $P_d(z,f) \equiv p_{l_1}(z) e^{\alpha_{l_1}(z)} + p_{l_2}(z) e^{\alpha_{l_2}(z)}$.
\end{theoF} 

\begin{rem}\label{r2.2}
First of all we would like to mention that Chen-Lian \cite{chen lian_BKMS20} did not point out the lacuna of {\em Theorem E}, but their result automatically rectifies {\em Theorem E}. Next we note that to prove Theorem F, the basic methods used by the authors are well known and is same as that adopted in 
several papers like {\em \cite{Latreuch_Mediterr J Math_2017,liao_ ann fenn13}}. The only innovative ideas exhibited by Chen-Lian \cite{chen lian_BKMS20} was to manipulate the notion of Cramer's rule in the proof. Unfortunately, at the time of execution of  Cramer's rule, Chen-Lian \cite{chen lian_BKMS20} made a mistake. Consider equation (3.2) {\em[}see P. 1067, \cite{chen lian_BKMS20}{\em ]}. 
In view of Crammer's rule, if $D_0\not= 0$, then only from equation (3.2), one can write equation (3.3) {\em[}see P. 1067, \cite{chen lian_BKMS20}{\em ]}. But in the line before (3.4), the authors used $D_0=0$ in (3.3) to deduce $D_1=0$, which is not at all acceptable and so to tackle this situation further investigations are needed. 
\end{rem}
Next let us accumulate the following points:- \newline 
i) From Remark 2.2 we see that the proof of {\em Theorem F} is incomplete and it is high time to dispel all the confusions cropped up from {\em Theorems E-F}.\newline  
ii) Considering all the results so far stated, a natural inquisition would be to investigate the case when the right hand side of the differential equation in {\it Theorem F} contains $k$-terms. 
\newline
iii) In 2011, Li \cite{ping li_jmaa11} proposed an open question that how to find the solutions of (\ref{e1.2}), where $p_1$, $p_2$ are constants and degree of the differential term $ P_d( z,f )$ is equal to $n -1$. But till now, without any extra supposition, nobody has been able to obtain any fruitful result in the literature. So Li's \cite{ping li_jmaa11} question is still open.\par
 In this respect, in view of previous points, considering the special case $f^{n-2}f' + P_d(z,f)$, it will be interesting to characterize the solutions of a more general form of differential equation, namely, \bea \label{e2.1} f^n+af^{n-2}f' + P_d(z,f) = \sum_{i=1}^{k}p_i(z)e^{\alpha_i(z)}, \eea where $a\in\mathbb{C}$, $P_d(z,f)$ is defined as in (\ref{e1.1}), $p_i(z)$ ($i = 1,2,\cdots,k$) are non-vanishing rational functions and $\alpha_i(z)$ ($i = 1,2,\cdots,k$) are distinct non-constant polynomials.\par
The above three points are the main motivations of writing this paper. Our result is the following: 

\begin{theo}\label{t2.1}  Consider the non-linear differential equation (\ref{e2.1}) with $\deg(\alpha_i-\alpha_j)\geq 1$ $(1\leq i\neq j\leq k)$.
	\begin{itemize}
		\item[(I)] Let $a=0$ and $d\leq n-k-1$, then the following two cases hold:
			\begin{itemize}
				\item [(IA)] Consider $k=1$ and $n\geq 2$. If (\ref{e2.1}) admits a meromorphic solution $f$ with finitely many poles, then $f$ is of the form
				$f (z) = q(z)e^{P(z)}$, where $q(z)$ is a non-vanishing rational function and $P(z)$ is a non-constant polynomial with $q^n(z) = p_1(z)$, $nP(z) = \alpha_1(z)$ and $P_d (z, f ) \equiv 0$.
				
				\item [(IB)] Consider $k\geq 2$ and $n \geq k+2$. If (\ref{e2.1}) has a meromorphic solution $f(z)$ with finitely many poles, then $\displaystyle\frac{\alpha'_i}{\alpha'_j}$ $(1\leq i\neq j\leq k)$ are rational numbers and $f(z)$ must be of the form:
				\beas f(z) = q(z)e^{P(z)},\eeas
				where $q(z)$ is a non-vanishing rational function and $P(z)$ is a non-constant
				polynomial.\par Also, if we rearrange $\{1,2,\ldots,k\}$ to $\{\tau_0,\tau_1,\ldots,\tau_{k-1}\}$ such that $p_i(z)e^{\alpha_i(z)}=p_{\tau_{i-1}}(z)e^{\alpha_{\tau_{i-1}}(z)}$ for $i=1,2,\ldots,k$, then $q^n=e^{-\mathcal{A}_0}p_{\tau_0},$ $\mathcal{A}_0$ is any constant with \beas nP'(z)=\alpha'_{\tau_0}(z)\;\;and\;\;P_d(z,f)=\sum\limits_{i=1}^{k-1}p_{\tau_i}(z)e^{\alpha_{\tau_i}(z)}.\eeas
			\end{itemize}
		\item[(II)] Next let $a\neq 0$ and $d\leq n-k-3$, then the following three cases hold:
			\begin{itemize}
				\item [(IIA)] Suppose $k=1$ and $n \geq 5$. Then the equation (\ref{e2.1}) does not admit any meromorphic solution with finitely many poles.
				
				\item [(IIB)] Suppose $k=2$ and $n \geq 6$. If (\ref{e2.1}) has a meromorphic solution $f(z)$ with finitely many poles, then $\displaystyle\frac{\alpha'_1}{\alpha'_2}$ are rational numbers and $f(z)$ must be of the form:
				\beas f(z) = q(z)e^{P(z)},\eeas
				where $q(z)$ is a non-vanishing rational function and $P(z)$ is a non-constant
				polynomial and $P_d(z,f)\equiv 0$. Also,
					\begin{itemize}
						\item [(i)] $q^n=e^{-\mathcal{B}_1}p_{1},\;aq^{n-2}(q'+qP')=e^{-\mathcal{B}_2}p_{2},\;\frac{\alpha'_{1}}{\alpha'_{2}}=\frac{n}{n-1},$ where $\mathcal{B}_i$ $(i=1,2)$ are constants
						\item [(ii)] or $q^n=e^{-\tilde{\mathcal{B}}_1}p_{1},\;aq^{n-2}(q'+qP')=e^{-\tilde{\mathcal{B}}_2}p_{2},\;\frac{\alpha'_{1}}{\alpha'_{2}}=\frac{n-1}{n},$ where $\tilde{\mathcal{B}}_i$ $(i=1,2)$ are constants.
					\end{itemize} 
				
				\item [(IIC)] Suppose $k\geq 3$ and $n \geq k+4$. If (\ref{e2.1}) has a meromorphic solution $f(z)$ with finitely many poles, then $\displaystyle\frac{\alpha'_i}{\alpha'_j}$ $(1\leq i\neq j\leq k)$ are rational numbers and $f(z)$ must be of the form:
				\beas f(z) = q(z)e^{P(z)},\eeas
				where $q(z)$ is a non-vanishing rational function and $P(z)$ is a non-constant
				polynomial.\par Also, if we rearrange $\{1,2,\ldots,k\}$ to $\{\mu,\nu,\kappa_1,\kappa_2, \ldots,\kappa_{k-2}\}$ such that $p_1(z)e^{\alpha_1(z)}=p_{\mu}(z)e^{\alpha_{\mu}(z)}$, $p_2(z)e^{\alpha_2(z)}=p_{\nu}(z)e^{\alpha_{\nu}(z)}$, $p_i(z)e^{\alpha_i(z)}=p_{\kappa_{i-2}}(z)e^{\alpha_{\kappa_{i-2}}(z)}$ for $i=3,4\ldots,k$, then $q^n=e^{-\mathcal{C}_{\mu}}p_{\mu}$ and $aq^{n-2}(q'+qP')=e^{-\mathcal{C}_{\nu}}p_{\nu}$, where $\mathcal{C}_{\mu},\;\mathcal{C}_{\nu}$ are any constants with \beas\hspace{5cc} nP'(z)=\alpha'_{\mu}(z),\;(n-1)P'(z)=\alpha'_{\nu}(z)\;and\;P_d(z,f)=\sum\limits_{i=1}^{k-2}p_{\kappa_i}(z)e^{\alpha_{\kappa_i}(z)}.\eeas
			\end{itemize}
	\end{itemize}	
\end{theo}

\begin{cor}
	Putting $k=3$, in {\em (IB)} of {\em Theorem \ref{t2.1}}, we have {\em Theorem F}. Therefore, our result is a huge improvement of {\em Theorem F}.
\end{cor}

The following examples show that all conclusions of {\it Theorem \ref{t2.1}} actually occurs for the cases $a=0$ and $a\not= 0$ respectively:

\begin{exm} Let $k=2$ and $n=4$, then $f(z)=\frac{z}{z+1}e^{z^2+2}$ satisfies the differential equation \beas f^4 + P_d(z,f) = e^3\left(\frac{z}{z+1}\right)^4e^{\alpha_1} + \frac{2z^2(z+1)+1}{(z+1)^2} e^{\alpha_2},\eeas where $P_d(z,f)=f'$, $\alpha_1=4z^2+5$ and $\alpha_2=z^2+2$. Clearly, $\frac{\alpha'_1}{\alpha'_2}$ is rational.
\end{exm}
\begin{exm} Let $k=2$ and $n=6$, then it is easy to verify that $f(z)= e^{\frac{2z}{3}}$ satisfies the differential equation \beas f^6 + af^4 f' = e^{4z} + \frac{2a}{3} e^{\frac{10z}{3}}.\eeas Note that here $P_d(z,f)\equiv 0$.
\end{exm}
\begin{exm} Let $k=3$ and $n=7$. Then choosing $P_d(z,f)=f''$, one can show that $f(z)= e^{\frac{2z}{7}}$ satisfies the differential equation \beas f^7 + af^5 f' + P_d(z,f) = e^{2z} + \frac{2a}{7} e^{\frac{12z}{7}} + \frac{4}{49} e^{\frac{2z}{7}}.\eeas \end{exm}

The next two examples show that in {\it (I)} of {\it Theorem \ref{t2.1}}, the bound  $d\leq n-k-1$ can not be extended to $d\leq n-k$ and it is the best possible estimation.
\begin{exm} Let $k=2, n=4$ and $d=2$, then it is easy to verify that $f(z)= e^z + z+1$ is a solution of the differential equation \beas f^4 + P_d(z,f) = e^{4z} + 4(z+1)e^{3z} ,\eeas where $P_d(z,f)=-6(z+1)^2(f'')^2 -3(z+1)^3 f'' -(z+1)^3 f$, but $ e^z + z+1$ can not be expressed as $q(z)e^{P(z)}$.
\end{exm} 
\begin{exm} Let $k=3, n=5$ and $d=2$, then it is easy to verify that $f(z)= e^z -1$ is a solution of the differential equation \beas f^5 + P_d(z,f) = e^{5z} -5 e^{4z} + 10e^{3z} ,\eeas where $P_d(z,f)=10 f f' + 5 f' + 1$, but $ e^z -1$ can not be written as $q(z)e^{P(z)}$.
\end{exm} 
Notice that, if $d\leq n-k-1$, then the bound $n\geq k+2$ can not be replaced by $n\geq k+1$, because in that case the equation (\ref{e2.1}) is not at all a differential equation. Now if we consider $n=k+1$ and $d= n-k$ then the next example shows that the conclusion (IB) of {\it Theorem \ref{t2.1}} cease to be hold. So we can say the bound for $n$ is the best possible. 
\begin{exm} Let $k=4, n=5$ and $d=1$. Here $f(z)= e^z +z-1$ is a solution of the differential equation \beas f^5 +P_d(z,f)= e^{5z} +5(z-1) e^{4z} + 10(z-1)^2e^{3z} + 10(z-1)^3 e^{2z} ,\eeas where $P_d(z,f)= - 4(z-1)^4 f'' - (z-1)^4 f$, but $ e^z +z-1$ can not be expressed as $q(z)e^{P(z)}$.
\end{exm}

The following examples show that in {\it (II)} of {\it Theorem \ref{t2.1}}, the bound  $d\leq n-k-3$ is sharp.

\begin{exm} Let $k=3, n=7$ and $d=2$, then it is easy to verify that $f(z)= e^z +1$ is a solution of the differential equation \beas f^7 - \frac{7}{2}f^5f' + P_d(z,f) = e^{7z} + \frac{7}{2}e^{6z} + \frac{7}{2} e^{5z},\eeas where $P_d(z,f)= - \frac{7}{2}ff' -1$, but $ e^z +1$ is not of the form $q(z)e^{P(z)}$. Here we note that $2=d>n-k-3=1$.
\end{exm}

\section{Lemmas.}
The following lemma can be easily derived from the proof of the Clunie lemma (see \cite{clunie_london math62,laine_book93}).
\begin{lem}\label{l3.1} Let $f(z)$ be a transcendental meromorphic solution of the differential equation \beas f^n(z) P(z,f) = Q(z,f),\eeas where $P(z,f), Q(z,f)$ are polynomials in $f$ and its derivatives such that the coefficients are small meromorphic functions of $f$, If the total degree of $Q(z, f)$ as a polynomial in $f$ and its derivatives is at most $n$, then \beas m\left(r,P(z,f)\right) = S(r,f), \eeas for all r out of a possible exceptional set of finite logarithmic measure. In particular, if $f$ is finite order then \beas m\left(r,P(z,f)\right) = O(\log\;r),\;\text{as}\; r\rightarrow \infty.\eeas 
\end{lem}

\begin{lem}\label{l3.2} \cite[Cramer's rule]{mirsky_book55} Consider the system of linear equation $AX= B$, where 	\beas A =\begin{pmatrix}[cccc] a_{11} & a_{12}& \ldots & a_{1n} \\  a_{21} & a_{22}& \ldots & a_{2n}\\  \vdots & \vdots &\ddots& \vdots\\ a_{n1} & a_{n2}& \ldots & a_{nn} \end{pmatrix},\; X= \begin{pmatrix}[c] x_1\\ x_2\\ \vdots\\ x_n \end{pmatrix}\;\text{and}\;B= \begin{pmatrix}[c] b_1\\ b_2\\ \vdots\\ b_n \end{pmatrix}.\eeas If $\det(A) \neq 0$, then the system has unique solution \beas \left(x_1,x_2,\cdots,x_n\right) = \left(\frac{\det(A_1)}{\det(A)}, \frac{\det(A_2)}{\det(A)}, \cdots, \frac{\det(A_n)}{\det(A)}\right),\eeas where \beas A_i =\begin{pmatrix}[cccccccc] a_{11} & a_{12}& \ldots &a_{1,i-1}& b_1 & a_{1,i+1} & \ldots & a_{1n} \\  a_{21} & a_{22}& \ldots &a_{2,i-1}& b_2 & a_{2,i+1} & \ldots & a_{2n}\\  \vdots & \vdots &\ddots& \vdots\\ a_{n1} & a_{n2}& \ldots &a_{n,i-1}& b_n & a_{n,i+1} & \ldots & a_{nn} \end{pmatrix}, \; \text{for}\; i=1,2,\ldots,n. \eeas\end{lem}

\begin{lem} \label{l3.3} \cite{yang yi_book03} Let $a_j(z)$ be entire function of finite order $\leq \rho$. Let $g_j(z)$ be entire and $g_k(z) - g_j(z)$, $(k\neq j)$ be a transcendental entire function or polynomial of degree greater than $\rho$. Then \beas \sum_{j=1}^{n}a_j(z)e^{g_j(z)} = a_0(z),\eeas holds only when \beas a_0(z) = a_1(z) = \cdots=a_n(z) \equiv 0.\eeas
\end{lem}
\begin{lem} \label{l3.4} \cite[Hadamard's factorization theorem]{yang yi_book03} Let $f(z)$ be a meromorphic function of finite order $\rho$ and \beas f(z) = a_kz^k + a_{k+1}z^{k+1} + \cdots \;\; (\text{where}\;a_k \neq 0) \eeas in the neighborhood of $z = 0$. Suppose that $b_1, b_2,\cdots$ are non-zero zeros of $f(z)$ and $c_1, c_2,\cdots$ are non-zero poles of $f(z)$, then \beas f(z) = z^ke^{Q(z)} \frac{P_1(z)}{P_2(z)},\eeas where $P_1(z)$ is a canonical product of non-zero zeros of $f(z)$, $P_2(z)$ is a canonical product of nonzero poles of $f(z)$, and $Q(z)$ is a polynomial of degree at most $\rho$.
\end{lem}


\section{Proof of the Theorem}

\begin{proof} [\bf\underline{Proof of Theorem \ref{t2.1}}] 
In view of {\em Remark \ref{r2.2}}, a detail proof is required for the sake of general reader.\par  To prove this theorem, we distinguish the following cases:\\
	\textbf{Case I:} Let $a=0$. Then we denote
	\bea\label{e4.1} h(z) = f^n + P_d(z, f).\eea
	Let $f$ be a rational solution of the non-linear differential equation (\ref{e2.1}). Then it is easy to see that	$h(z)$ is a small function of $p_i(z)e^{\alpha_i(z)}$ $(i = 1, 2, \ldots,k)$. As $\deg\{\alpha_i(z)-\alpha_j(z)\}\geq 1$, so by {\it{Lemma \ref{l3.3}}}, we get $p_i(z)\equiv 0$ $(i = 1, 2, \ldots, k)$, which contradicts $p_i(z)$ $(i = 1, 2, \ldots, k)$ are non-vanishing rational functions. Hence $f$ must be a transcendental.\par
	Now we will show that order of $f$ is finite. As $f$ has finitely many poles, so \beas nT(r,f)&=&m(r,f^n)+S(r,f)\\&\leq& T\left(r,\sum\limits_{i=1}^{k}p_i(z)e^{\alpha_i(z)}\right)+m(r, P_d(z, f))+S(r,f)\\&\leq&Ar^{\eta}+dT(r,f)+S(r,f), \eeas where $\eta=\max\{\deg{\alpha_1},\deg{\alpha_2},\ldots,\deg{\alpha_k}\}.$ Hence $(n-d)T(r,f)\leq Ar^{\eta}+S(r,f)$ and $f$ is of finite order.\par
	\textbf{Sub-case IA:} First suppose that $k=1$, then the result follows from {\it Theorem 1.7} of \cite{liao_ compl ver elliptic 15}.\par
	\textbf{Sub-case IB:} Next suppose that $k\geq 2$. For $k=2$, we refer {\it Theorem D}. So we assume $k\geq 3$.
	 Differentiating (\ref{e2.1}), $k-1$ times, we get the following system of equations
	\bea\label{e4.2}\left.\begin{array}{clcr}  h&=&\sum\limits_{i=1}^{k}p_ie^{\alpha_i},  \\   h'&=&\sum\limits_{i=1}^{k}\left[p_i'+p_i\alpha_i'\right]e^{\alpha_i}, \\ h''&=&\sum\limits_{i=1}^{k}\left[p_i''+2p_i'\alpha_i'+p_i\alpha_i''+p_i(\alpha_i')^2\right]e^{\alpha_i},\\ &\vdots& \\ h^{(k-1)}&=&\sum\limits_{i=1}^{k}\left[p_i^{(k-1)}+p_i^{(k-2)}Q_{1}(\alpha_i')+p_i^{(k-3)}Q_{2}(\alpha_i',\alpha_i'')\right.\\&&\left.+\ldots+p_iQ_{k-1}(\alpha_i',\alpha_i'',\ldots,\alpha_i^{(k-1)})\right]e^{\alpha_i},\end{array}\right\}\eea 
	where $Q_{j}(\alpha_i',\alpha_i'',\ldots,\alpha_i^{(j)})$ are differential polynomials of $\alpha_i'$ with degree $j$ $(j=1,2,\ldots,k-1;\;i=1,2,\ldots,k)$.
	From the system (\ref{e4.2}), the determinant of the coefficient matrix is
	
	\beas D_0=\begin{vmatrix} p_1 & \ldots & p_k \\p_1'+p_1\alpha_1' &\ldots& p_k'+p_k\alpha_k' \\  \vdots &\ddots& \vdots\\ \begin{matrix} p_1^{(k-1)}+p_1^{(k-2)}Q_{1}(\alpha_1')+\ldots\\+ p_1Q_{k-1}(\alpha_1',\alpha_1'',\ldots,\alpha_1^{(k-1)})\end{matrix} &\ldots& \begin{matrix}	p_k^{(k-1)}+p_k^{(k-2)}Q_{1}(\alpha_k')+\ldots\\+ p_kQ_{k-1,k-1}(\alpha_k',\alpha_k'',\ldots,\alpha_k^{(k-1)})\end{matrix}\end{vmatrix}.\eeas It is clear that $D_0$ is a rational function. Let us consider 
	
	\beas D_1=\begin{vmatrix} h & p_2 & \ldots & p_k \\ h' & p_2'+p_2\alpha_2' &\ldots& p_k'+p_k\alpha_k' \\  \vdots & \vdots &\ddots& \vdots\\    h^{(k-1)} & \begin{matrix}	p_2^{(k-1)}+p_2^{(k-2)}Q_{1}(\alpha_2')+\ldots+\\ p_2Q_{k-1}(\alpha_2',\alpha_2'',\ldots,\alpha_2^{(k-1)})\end{matrix} &\ldots& \begin{matrix}	p_k^{(k-1)}+p_k^{(k-2)}Q_{1}(\alpha_k')+\ldots+\\ p_kQ_{k-1}(\alpha_k',\alpha_k'',\ldots,\alpha_k^{(k-1)})\end{matrix}\end{vmatrix}.\eeas
	 Now we distinguish the following two cases:\\
	 	
	\textbf{Sub-case IB.1:} If $D_0\equiv 0$, then the rank of the coefficient matrix of the system (\ref{e4.2}) is equal to $m\leq k-1$. As the system of equation has a solution, so the rank of the augmented matrix,
	
	\beas \tilde{D_0}=\begin{pmatrix}[ccc|c] p_1 & \ldots & p_k & h \\  p_1'+p_1\alpha_1' &\ldots& p_k'+p_k\alpha_k' &h'\\  \vdots & \ddots &\vdots& \vdots\\  \begin{matrix} p_1^{(k-1)}+p_1^{(k-2)}Q_{1}(\alpha_1')+\ldots+\\ p_1Q_{k-1}(\alpha_1',\alpha_1'',\ldots,\alpha_1^{(k-1)})\end{matrix} &\ldots& \begin{matrix}	p_k^{(k-1)}+p_k^{(k-2)}Q_{1}(\alpha_k')+\ldots+\\ p_kQ_{k-1}(\alpha_k',\alpha_k'',\ldots,\alpha_k^{(k-1)})\end{matrix}  & h^{(k-1)} \end{pmatrix},\eeas must be equal to $m\leq k-1$. So, all $k\times k$ minors of $\tilde{D_0}$ are zero, which implies 
	
	\beas \begin{vmatrix} p_2 & \ldots & p_k & h \\  p_2'+p_2\alpha_2' &\ldots& p_k'+p_k\alpha_k' &h'\\  \vdots & \ddots &\vdots& \vdots\\ \begin{matrix}	p_2^{(k-1)}+p_2^{(k-2)}Q_{1}(\alpha_2')+\ldots\\+ p_2Q_{k-1}(\alpha_2',\alpha_2'',\ldots,\alpha_2^{(k-1)})\end{matrix} &\ldots& \begin{matrix}	p_k^{(k-1)}+p_k^{(k-2)}Q_{1}(\alpha_k')+\ldots\\+ p_kQ_{k-1}(\alpha_k',\alpha_k'',\ldots,\alpha_k^{(k-1)})\end{matrix} & h^{(k-1)} \end{vmatrix}=0,\eeas which implies, $D_1\equiv 0$.\\
	So, \bea\label{e4.3} M_{11}h-M_{21}h'+\ldots+(-1)^{k-1}M_{k1}h^{(k-1)}=0,\eea
	where \beas M_{11}=\begin{vmatrix} p_2'+p_2\alpha_2' &\ldots& p_k'+p_k\alpha_k' \\   \vdots &\ddots& \vdots\\   \begin{matrix}	p_2^{(k-1)}+p_2^{(k-2)}Q_{1}(\alpha_2')+\ldots\\+ p_2Q_{k-1}(\alpha_2',\alpha_2'',\ldots,\alpha_2^{(k-1)})\end{matrix} &\ldots& \begin{matrix}	p_k^{(k-1)}+p_k^{(k-2)}Q_{1}(\alpha_k')+\ldots\\+ p_kQ_{k-1}(\alpha_k',\alpha_k'',\ldots,\alpha_k^{(k-1)})\end{matrix}\end{vmatrix},\eeas
	
	\beas M_{21}= \begin{vmatrix} p_2 & \ldots & p_k \\ p_2''+2p_2'\alpha_2'+p_2\alpha_2''+p_2(\alpha_2')^2 &\ldots &p_k''+2p_k'\alpha_k'+p_k\alpha_k''+p_k(\alpha_k')^2 \\ \vdots &\ddots& \vdots\\  \begin{matrix}	p_2^{(k-1)}+p_2^{(k-2)}Q_{1}(\alpha_2')+\ldots\\+ p_2Q_{k-1}(\alpha_2',\alpha_2'',\ldots,\alpha_2^{(k-1)})\end{matrix} &\ldots& \begin{matrix}	p_k^{(k-1)}+p_k^{(k-2)}Q_{1}(\alpha_k')+\ldots\\+ p_kQ_{k-1}(\alpha_k',\alpha_k'',\ldots,\alpha_k^{(k-1)})\end{matrix}\end{vmatrix}\eeas
	and so 
	
	\beas M_{k1}=\begin{vmatrix} p_2 & \ldots & p_k \\  p_2'+p_2\alpha_2' &\ldots& p_k'+p_k\alpha_k' \\ \vdots &\ddots& \vdots\\  \begin{matrix}	p_2^{(k-2)}+p_2^{(k-3)}Q_{1}(\alpha_2')+\ldots\\+ p_2Q_{k-2}(\alpha_2',\alpha_2'',\ldots,\alpha_2^{(k-2)})\end{matrix} &\ldots& \begin{matrix}	p_k^{(k-2)}+p_k^{(k-3)}Q_{1}(\alpha_k')+\ldots\\+ p_kQ_{k-2}(\alpha_k',\alpha_k'',\ldots,\alpha_k^{(k-2)})\end{matrix}\end{vmatrix}\eeas
	
	are rational functions.\\
	Substituting the expression (\ref{e4.1}) of $h(z)$ into (\ref{e4.3}), we get
	\bea\label{e4.4} M_{11}f^n-M_{21}(f^n)'+\ldots+(-1)^{k-1}M_{k1}(f^n)^{(k-1)}=Q_1,\eea
	where $$Q_1=-\left[M_{11}P_d(z,f)-M_{21}P_d^{\prime}(z,f)+\ldots+(-1)^{k-1}M_{k1}P_d^{(k-1)}(z,f)\right]$$ is a differential polynomial in $f$ with rational functions as its coefficients and degree of $Q_1 \leq d$.\\
	One can easily check that \bea\label{e4.5}(f^n)^{(t)}&=&(nf^{n-1}f')^{(t-1)}=n\sum\limits_{i=0}^{t-1}\binom{t-1}{i}(f^{n-1})^{(i)}f^{(t-i)}\nonumber\\&=&n\left[f^{n-1}f^{(t)}+(t-1)(n-1)f^{n-2}f'f^{(t-1)}\right.\nonumber\\&&\left.+\sum\limits_{i=2}^{t-1}\binom{t-1}{i}f^{(t-i)}\left\{(n-1)f^{n-2}f^{(i)}+\sum\limits_{j=2}^{i-1}\sum\limits_{\lambda}\gamma_{j\lambda}f^{n-j-1}(f')^{\lambda_{j,1}}\right.\right.\nonumber\\&&\left.\left.(f'')^{\lambda_{j,2}}\ldots(f^{(i-1)})^{\lambda_{j,i-1}}+(n-1)(n-2)\ldots(n-i)f^{n-i-1}(f')^{i} \right\}\right],\eea
	where $\gamma_{j\lambda}$ are positive integers, $\lambda_{j,1},\ldots,\lambda_{j,i-1}$ are non-negative integers and sum $\sum\limits_{\lambda}$ is carried out such that \bea\label{e4.6}\lambda_{j,1}+\lambda_{j,2}+\ldots+\lambda_{j,i-1}=j\;\;and\;\;\lambda_{j,1}+2\lambda_{j,2}+\ldots+(i-1)\lambda_{j,i-1}=i.\eea Now we define  \bea\label{e4.7}\xi_t&=&\frac{(f^n)^{(t)}}{f^{n-k+1}},\eea for $t=1,2,\ldots,k-1$; $k\geq 3$.
	Using (\ref{e4.7}) in (\ref{e4.4}), we have 
	
	\bea\label{e4.8} f^{n-k+1}R_1=Q_1,\eea where \bea\label{e4.9}R_1=M_{11}f^{k-1}-M_{21}\xi_1+\ldots+(-1)^{k-1}M_{k1}\xi_{k-1}.\eea
	As $d \leq n-k-1$ and $f$ is finite order, so combining (\ref{e4.8}) with {\it{Lemma \ref{l3.1}}} we get
	$$m(r,R_1) = O(\log r).$$ On the other hand, as we assume, $f$ has finitely many poles, so we have \beas T(r,R_1) = m(r,R_1) + N(r,R_1) = O(\log r),\eeas i.e., $R_1$ is a rational function.	Next we study the following two subcases:\\
		\textbf{Sub-case IB.1.1:} If $R_1(z)\equiv 0$, then from (\ref{e4.9}) we have
	\bea\label{e4.10} M_{11}f^{k-1}=-\left[-M_{21}\xi_1+\ldots+(-1)^{k-1}M_{k1}\xi_{k-1}\right].\eea
	
	We will show that $f$ has at most finitely many zeros. On the contrary, suppose that $f$ has infinitely many zeros. So we can consider a point $z_1$ such that
     $f(z_1) = 0$ but $z_1$ is neither a zero nor a pole of $M_{j1}$ $(j=1,2,\ldots,k)$. \par
	Now let, $M_{k1}\not\equiv 0$. From the construction of $\xi_{k-1}$ and (\ref{e4.5}) we know $\xi_{k-1}$ contains one term, corresponding to $i=k-2$ in which, power of $f$ is $0$. Considering this term in $\xi_{k-1}$, we see that $f'(z_1)=0$, which implies $z_1$ is a multiple zero of $f$ of multiplicity $p_1\geq2$. It follows that $z_1$ is a zero of the left hand side of (\ref{e4.10}) with multiplicity $(k-1)p_1$, where as it is a zero of the right hand side of (\ref{e4.10}) with multiplicity $(k-1)(p_1-1)$, a contradiction.\par
	Next, assume $M_{k1}\equiv 0$. If $z_1$ is a simple zero, then $z_1$ is a zero with multiplicity $k-1$ of left hand side of (\ref{e4.10}) and a zero with multiplicity $1$ of right hand side of (\ref{e4.10}), a contradiction. If $z_1$ is a multiple zero with multiplicity $q_1\geq2$, then $z_1$ is a zero with multiplicity $(k-1)q_1$ of left hand side of (\ref{e4.10}) and a zero with multiplicity $(k-1)q_1-(k-2)$ of right hand side of (\ref{e4.10}), a contradiction. Thus, $f$ has at most finitely many zeros.\\
		\textbf{Sub-case IB.1.2:} If $R_1(z)\not\equiv 0$, then (\ref{e4.8}) becomes \bea\label{e4.11}f^{n-k}(fR_1)=Q_1.\eea By {\it{Lemma \ref{l3.1}}}, we have
	$$m(r,fR_1) = O(\log r).$$
	From our assumption, since $fR_1$ has finitely many poles. Then
	$$T(r,fR_1) = m(r,fR_1) + N(r,fR_1) = O(\log r).$$
	Therefore, $fR_1$ is a rational function, which contradicts that $f$ is transcendental.\\	
	\textbf{Sub-case IB.2:} If $D_0\not\equiv 0$, then by  {\it {Lemma \ref{l3.2}}} we get 
	\bea\label{e4.12} D_0 e^{\alpha_1} = D_1. \eea Differentiating (\ref{e4.12}), we have \bea\label{e4.13}(D_0'+D_0\alpha_1')e^{\alpha_1} = D_1'.\eea
	Eliminating $e^{\alpha_1}$ from (\ref{e4.12}) and (\ref{e4.13}), we get
	\bea\label{e4.14} D_1'D_0-D_1D_0'=\alpha_1'D_1D_0.\eea
	Substituting $D_1=M_{11}h-M_{21}h'+\ldots+(-1)^{k-1}M_{k1}h^{(k-1)}$ and $D_1'=M'_{11}h+(M_{11}-M'_{21})h'-(M_{21}-M'_{31})h''+\ldots+(-1)^{k-2}(M_{k-1\; 1}-M'_{k 1})h^{(k-1)}+(-1)^{k-1}M_{k 1}h^{(k)}$ in (\ref{e4.14}), we have
	\bea\label{e4.15} A_1h+A_2h'+\ldots+A_{k+1}h^{(k)}=0,\eea 
	where \beas\begin{array}{clcr} A_1&=&M'_{11}D_0-M_{11}(D_0'+\alpha_1'D_0),\hspace{9.7cc}\\ A_2&=&(M_{11}-M_{21}')D_0+M_{21}(D_0'+\alpha_1'D_0),\hspace{6.8cc} \\ &\vdots& \\ A_k&=&(-1)^{k-2}(M_{k-1\;1}-M'_{k\;1})D_0-(-1)^{k-1}M_{k1}(D_0'+\alpha_1'D_0),\\A_{k+1}&=&(-1)^{k-1}M_{k1}D_0\hspace{13.8cc}\end{array}\eeas 
	are rational functions.\par Substituting the expressions of $h',h'',\ldots,h^{(k)}$ into (\ref{e4.15}), we get
	\bea\label{e4.16} A_1f^n+A_2(f^n)'+\ldots+A_k(f^n)^{(k-1)}+A_{k+1}(f^n)^{(k)}=Q_2,\eea
	where $Q_2 = -[A_1P_d(z,f)+A_2P_d'(z,f)+\ldots+A_kP_d^{(k-1)}(z,f)+A_{k+1}P_d^{(k)}(z,f)]$ is a differential polynomial in $f$ with rational functions as its coefficients and degree of $Q_2 \leq d$. \par
	Using (\ref{e4.5}) we define \bea\label{e4.17} \psi_t=\frac{(f^n)^{(t)}}{f^{n-k}},\eea $t=1,2,\ldots,k$; $k\geq 3$. Now applying (\ref{e4.5}) and (\ref{e4.17}) in (\ref{e4.16}), we have \bea\label{e4.18} f^{n-k}R_2=Q_2,\eea where \bea\label{e4.19}R_2=f^kA_1+A_2\psi_1+\ldots+A_{k+1}\psi_k.\eea\\
	Noting the fact $d \leq n-k-1$ and combining (\ref{e4.18}) with {\it{Lemma \ref{l3.1}}} we obtain
	$$m(r,R_2) = O(\log r).$$ By the assumption, $f$ has finitely many poles, then $$T(r,R_2) = m(r,R_2) + N(r,R_2) = O(\log r),$$ i.e., $R_2$ is a rational function.\\
	Next we discuss two sub cases as follows:\\	
	\textbf{Sub-case IB.2.1:} If $R_2(z)\equiv 0$, then by (\ref{e4.19}) we have
	\bea\label{e4.20} f^kA_1=-(A_2\psi_1+\ldots+A_{k+1}\psi_k).\eea
	Now proceeding in the same way as done in {\it Sub-case IB.1.1} and replacing $k-1$ by $k$ we can show that $f$ has finitely many zeros. So we omit the details.
	
	\textbf{Sub-case IB.2.2:} If $R_2(z)\not\equiv 0$, then (\ref{e4.18}) becomes \bea\label{e4.21}f^{n-k-1}(fR_2)=Q_2.\eea Next similar to subcase {\it Sub-case IB.1.2} we can get a contradiction.
	
	So from the above discussion we conclude that $f$ is a transcendental meromorphic function with finite zeros and poles. Now by {\it{Lemma \ref{l3.4}}}, we can say that
	\bea\label{e4.22}f(z) = q(z)e^{P(z)},\eea where $q(z)$ is a non-vanishing rational function, and $P(z)$ is a non-constant polynomial. Substituting (\ref{e4.22}) into (\ref{e2.1}) yields
	\bea\label{e4.23} q^n(z)e^{nP(z)}+\sum\limits_{j=0}^{d}\beta_j(z)e^{jP(z)}=\sum\limits_{i=1}^{k}p_i(z)e^{\alpha_i(z)},\eea
	where $\beta_j(z)$ are rational functions.\\
	
	 Now we can rearrange $\{1,2,\ldots,d\}$ to $\{\sigma_1,\sigma_2,\ldots,\sigma_d\}$ such that $\beta_j(z)e^{jP(z)}=\beta_{\sigma_j}(z)e^{{\sigma_j}P(z)}$ for $j=1,\ldots,d$ and $\{1,2,\ldots,k\}$ to $\{\tau_0,\tau_1,\ldots,\tau_{k-1}\}$ such that $p_i(z)e^{\alpha_i(z)}=p_{\tau_{i-1}}(z)e^{\alpha_{\tau_{i-1}}(z)}$ for $i=1,2,\ldots,k$. Then (\ref{e4.23}) can be written as \bea\label{e4.24} q^n(z)e^{nP(z)}+\sum\limits_{j=0}^{d}\beta_{\sigma_j}(z)e^{{\sigma_j}P(z)}=\sum\limits_{i=1}^{k}p_{\tau_{i-1}}(z)e^{\alpha_{\tau_{i-1}}(z)}.\eea 
	
	Since $n>d$, $\deg(\alpha_{\tau_i}-\alpha_{\tau_j})\geq 1$ $(1\leq \tau_i\neq \tau_j\leq k)$ and $q$, $p_i$ $(i=1,2,\ldots,k)$ are all non-zero rational functions, then using {\it {Lemma \ref{l3.3}}} on (\ref{e4.24}) we have
	 \beas nP=\alpha_{\tau_0}+\mathcal{A}_0,\sigma_1P=\alpha_{\tau_1}+\mathcal{A}_1,\sigma_2P=\alpha_{\tau_2}+\mathcal{A}_2,\ldots,\sigma_{k-1}P=\alpha_{\tau_{k-1}}+\mathcal{A}_{k-1},\eeas where $\mathcal{A}_i$ $(i=0,1,\ldots,k-1)$ are constants such that \beas q^n=e^{-\mathcal{A}_0}p_{\tau_0},\beta_{\sigma_1}=e^{-\mathcal{A}_1}p_{\tau_1},\beta_{\sigma_2}=e^{-\mathcal{A}_2}p_{\tau_2},\ldots,\beta_{\sigma_{k-1}}=e^{-\mathcal{A}_{k-1}}p_{\tau_{k-1}}.\eeas
	Also, $\beta_0(z)=0$ and $\beta_{\sigma_j}(z)\equiv 0$ for all $\sigma_j\neq \sigma_1,\sigma_2,\ldots,\sigma_{k-1}$ with $0\leq \sigma_j\leq d$.\par Therefore, \beas \alpha'_{\tau_0}:\alpha'_{\tau_1}:\ldots:\alpha'_{\tau_{k-1}}=n:\sigma_1:\ldots:\sigma_{k-1},\eeas
	\beas nP'=\alpha'_{\tau_0}\;\;and\;\;P_d(z,f)=\sum\limits_{i=1}^{k-1}p_{\tau_i}(z)e^{\alpha_{\tau_i}(z)}.\eeas
	
	\textbf{Case II:} Let $a\neq 0$. In this case we consider \bea\label{e4.25}\tilde{h}(z)=f^n+af^{n-2}f' + P_d(z,f)\eea and similar as {\it Case I}, we can prove that $f$ is a finite order transcendental meromorphic function.\par
	
	\textbf{Sub-case IIA:} Suppose that $k=1$. Then (\ref{e2.1}) becomes \bea \label{e4.26} f^n+af^{n-2}f' + P_d = p_1e^{\alpha_1}.\eea Differentiating (\ref{e4.26}) we have, \bea \label{e4.27} nf^{n-1}f'+a(n-2)f^{n-3}(f')^2+af^{n-2}f''+ P'_d = p_1\left(\frac{p'_1}{p_1}+\alpha'_1\right)e^{\alpha_1}.\eea Eliminating $e^{\alpha_1}$ from (\ref{e4.26}) and (\ref{e4.27}), we have \bea\label{e4.28} f^{n-3}R_3=Q_3,\eea
	where \bea\label{e4.29}R_3=nf^2f'+a(n-2)(f')^2+aff''-\left(\alpha'_1+\frac{p'_1}{p_1}\right)(f^3+aff')\eea and $Q_3=\left(\alpha'_1+\frac{p'_1}{p_1}\right)P_d-P'_d$ with degree $d$. 
	Since $d\leq n-k-3 = n-4$, so from (\ref{e4.29}) and {\it{Lemma \ref{l3.1}}}, we have $$m(r,R_3)=O(\log r).$$ By the hypothesis, $f$ has finitely many poles, thus
	\beas T(r,R_3) = m(r,R_3) + N(r,R_3) = O(\log r),\eeas i.e., $R_3$ is a rational function.\par
	\textbf{Sub-case IIA.1:} Let $R_3\equiv 0$. Then from (\ref{e4.29}) we have \bea\label{e4.30} \left(\alpha'_1+\frac{p'_1}{p_1}\right)f^3=nf^2f'+a(n-2)(f')^2+aff''-a\left(\alpha'_1+\frac{p'_1}{p_1}\right)ff'.\eea
	\par Now we will prove that $f$ has only finitely many zeros. On the contrary, suppose that $f$ has infinitely many zeros. It is clear from (\ref{e4.30}) that all zeros of $f$ are multiple zeros. Assume $z_3$ be a zero of $f$ but not the zeros or poles of the coefficients of (\ref{e4.30}), with multiplicity $p_3\geq 2$. Then comparing the multiplicity of $f$ on both side of (\ref{e4.30}) we have, $z_3$ is a zero of the left hand side of (\ref{e4.30}) with multiplicity $3p_3$ and a zero of the right hand side of (\ref{e4.30}) with multiplicity $2(p_3-1)$, which is a contradiction. So, $f$ has at most finitely many zeros.\par Hence $f$ is a transcendental meromorphic function with finite zeros and poles. Now by {\it{Lemma \ref{l3.4}}}, we can say that $f(z) = q(z)e^{P(z)}$, where $q(z)$ is a non-vanishing rational function and $P(z)$ is a non-constant polynomial. Substituting the form of $f$ into (\ref{e4.26}) yields
	\bea\label{e4.31} q^n(z)e^{nP(z)}+aq^{n-2}(q'+qP')e^{(n-1)P(z)}+\sum\limits_{j=0}^{d}\eta_j(z)e^{jP(z)}=p_1(z)e^{\alpha_1(z)},\eea
	where $\eta_j(z)$ are rational functions.  As $q \neq 0$, so by applying {\it Lemma \ref{l3.3}} on (\ref{e4.31}) we get $q'+qP' = 0$, i.e., $q(z)=D_1/e^{P(z)}$ for a constant $D_1$, which is a contradiction that $q$ is rational.\par
	\textbf{Sub-case IIA.2:} Let $R_3\not\equiv 0$. Then from (\ref{e4.28}) we have, \bea\label{e4.32} f^{n-4}(fR_3)=Q_3.\eea
	Since $d\leq n-4$, combining (\ref{e4.32}) with {\it{Lemma \ref{l3.1}}}, we have $$m(r,fR_3)=O(\log r).$$ On the other hand, as we assume $f$ has finitely many poles, thus
	\beas T(r,fR_3) = m(r,fR_3) + N(r,fR_3) = O(\log r),\eeas i.e., $fR_3$ is a rational function and we have $R_3$ is rational, but $f$ is transcendental, a contradiction.\\ Therefore, in this case, there does not exist any transcendental meromorphic solution.\par 
	
	\textbf{Sub-case IIB:} In this case $k=2$. Then (\ref{e2.1}) becomes \bea \label{e4.33} f^n+af^{n-2}f' + P_d = p_1e^{\alpha_1}+p_2e^{\alpha_2}.\eea 
	Differentiating (\ref{e4.33}) and eliminating $e^{\alpha_2}$ we get 
	\bea\label{e4.34}&&\left(\frac{p_2'}{p_2}+\alpha'_2\right)f^n+a\left(\frac{p'_2}{p_2}+\alpha'_2\right)f^{n-2}f'-nf^{n-1}f'-a(n-2)f^{n-3}(f')^2-af^{n-2}f''\nonumber\\&&-P'_d+\left(\frac{p_2'}{p_2}+\alpha'_2\right)P_d=p_1Ae^{\alpha_1},\eea where $A\equiv \left(\frac{p'_2}{p_2}-\frac{p'_1}{p_1}+\alpha'_2-\alpha'_1\right)\not\equiv 0$, otherwise it contradicts our assumption $\deg(\alpha_i-\alpha_j)\geq 1$ $(1\leq i\neq j\leq k)$.\\
	Now differentiating (\ref{e4.34}) and eliminating $e^{\alpha_1}$, we have
	\bea\label{e4.35} f^{n-4}R_4=Q_4,\eea
	where \bea\label{e4.36}  R_4&=& h_1(z)(f^4+af^{2}f')+h_2(z)(nf^{3}f'+af^{2}f''+a(n-2)f(f')^2)+af^{2}f^{(3)}\nonumber\\&&+3a(n-2)ff'f''+n(n- 1)f^{2}(f')^2+ nf^{3}f''+ a(n - 2)(n -3)(f')^3 \eea and  $$Q_4=-P''_d -h_2(z)P'_d-h_1(z)P_d$$ such that $$h_1(z)=\left(\frac{p_2'}{p_2}+\alpha'_2\right)\left(\frac{p_1'}{p_1}+\frac{A'}{A}+\alpha'_1-\frac{\left(\frac{p_2'}{p_2}+\alpha'_2\right)'}{\left(\frac{p_2'}{p_2}+\alpha'_2\right)}\right),$$ $$h_2(z)=-\left(\frac{p_1'}{p_1}+\frac{A'}{A}+\alpha'_1+\frac{p_2'}{p_2}+\alpha'_2\right).$$
	Then it follows from {\it Lemma \ref{l3.1}} that $R_4$ is a rational function. Now we consider the following two cases:\par
	\textbf{Sub-case IIB.1:} Let $R_4=0$.\\ Then (\ref{e4.36}) can be written as \bea\label{e4.37} h_1(z)(f^4+af^{2}f')&=&-[h_2(z)(nf^{3}f'+af^{2}f''+a(n-2)f(f')^2)+af^{2}f^{(3)}+3a(n-2)ff'f''\nonumber\\&&+n(n- 1)f^{2}(f')^2+ nf^{3}f''+ a(n - 2)(n -3)(f')^3] \eea
	\par If $f$ has infinitely many zeros, it follows from (\ref{e4.37}) that zeros of $f$ are of
	multiplicity $p_4\geq2$. Let $z_4$ be a zero of $f$ but not zeros or poles of the coefficients. Then comparing the multiplicity of $f$ on both side of (\ref{e4.37}) we have, $z_4$ is a zero of the left hand side of (\ref{e4.37}) with multiplicity $3p_4-1$ and a zero of the right hand side of (\ref{e4.37}) with multiplicity $3(p_4-1)$, which is a contradiction. So, $f$ has at most finitely many zeros.\par So, $f$ is a transcendental meromorphic function with finite zeros and poles. Now by {\it{Lemma \ref{l3.4}}}, we can say that $f(z) = q(z)e^{P(z)}$, where $q(z)$ is a non-vanishing rational function and $P(z)$ is a non-constant polynomial. Substituting the form of $f$ into (\ref{e4.33}) yields
	\bea\label{e4.38} q^n(z)e^{nP(z)}+aq^{n-2}(q'+qP')e^{(n-1)P(z)}+\sum\limits_{j=0}^{d}\zeta_j(z)e^{jP(z)}=p_1(z)e^{\alpha_1(z)}+p_2(z)e^{\alpha_2(z)},\eea
	where $\zeta_j(z)$ are rational functions. Applying {\it Lemma \ref{l3.3}} on (\ref{e4.38}), we get \beas nP=\alpha_{1}+\mathcal{B}_1,\;\;(n-1)P=\alpha_{2}+\mathcal{B}_2,\eeas where $\mathcal{B}_i$ $(i=1,2)$ are constants such that \beas q^n=e^{-\mathcal{B}_1}p_{1},\;\;aq^{n-2}(q'+qP')=e^{-\mathcal{B}_2}p_{2}\eeas or \beas nP=\alpha_{2}+\tilde{\mathcal{B}}_1,\;\;(n-1)P=\alpha_{1}+\tilde{\mathcal{B}}_2,\eeas where $\tilde{\mathcal{B}}_i$ $(i=1,2)$ are constants such that \beas q^n=e^{-\tilde{\mathcal{B}}_1}p_{1},\;\;aq^{n-2}(q'+qP')=e^{-\tilde{\mathcal{B}}_2}p_{2}.\eeas Also, $\zeta_j(z)=0$ for all $j=0,1,\ldots,d$, i.e., $P_d(z,f)\equiv 0$.	Clearly, \beas \frac{\alpha'_{1}}{\alpha'_{2}}=\frac{n}{n-1}\;\;and\;\;\frac{\alpha'_{1}}{\alpha'_{2}}=\frac{n-1}{n}.\eeas
	\textbf{Sub-case IIB.2:} Let $R_4\neq 0$.\\
	Then from (\ref{e4.35}) we have, \bea\label{e4.39} f^{n-5}(fR_4)=Q_4.\eea
	Since $d\leq n-5$, combining (\ref{e4.39}) with {\it{Lemma \ref{l3.1}}}, we have $$m(r,fR_4)=O(\log r).$$ By the hypothesis, $f$ has finitely many poles, thus
	\beas T(r,fR_4) = m(r,fR_4) + N(r,fR_4) = O(\log r),\eeas i.e., $fR_4$ is a rational function and we have $R_4$ is rational, contradiction that $f$ is transcendental.\\	
	\textbf{Sub-case IIC:} In this case, we suppose that $k\geq 3$. We consider $\tilde{h}(z)$ in (\ref{e4.25}) instead of $h(z)$ in (\ref{e4.1}) and proceed similarly as {\it Sub-case IB} with necessary changes.\\
	\textbf{Sub-case IIC.1:} Now similar as {\it Sub-case IB.1} we proceed upto (\ref{e4.4}). In this case, (\ref{e4.4}) becomes \bea\label{e4.40} M_{11}(f^n+af^{n-2}f')+ \cdots+(-1)^{k-1}M_{k1}(f^n+af^{n-2}f')^{(k-1)}=Q_1\eea and similar as (\ref{e4.5}), we can check that \bea\label{e4.41}(f^{n-2}f')^{(t)}&=&f^{n-2}f^{(t+1)}+t(n-2)f^{n-3}f'f^{(t)}\nonumber\\&&+\sum\limits_{i=2}^{t}\binom{t}{i}f^{(t-i+1)}\left\{(n-2)f^{n-3}f^{(i)}+\sum\limits_{j=2}^{i-1}\sum\limits_{\lambda}\gamma_{j\lambda}f^{n-j-2}(f')^{\lambda_{j,1}}\right.\nonumber\\&&\left.(f'')^{\lambda_{j,2}}\ldots(f^{(i-1)})^{\lambda_{j,i-1}}+(n-2)(n-3)\ldots(n-i-1)f^{n-i-2}(f')^{i}\right\}.\hspace{2cc}\eea
	In this case, using (\ref{e4.5}) and (\ref{e4.41}), we define  \bea\label{e4.42}\xi_t&=&\frac{(f^n+af^{n-2}f')^{(t)}}{f^{n-k-1}},\eea for $t=1,2,\ldots,k-1$; $k\geq 3$. Then using (\ref{e4.42}), (\ref{e4.8}) changes to \bea\label{e4.43} f^{n-k-1}R_5=Q_5.\eea After that proceeding similarly, we can prove that $R_5$ is rational and can consider the following two subcases.\\	
	\textbf{Sub-case IIC.1.1:} If $R_5(z)\equiv 0$, then (\ref{e4.9}) changes to
	\bea\label{e4.44} M_{11}(f^{k+1}+af^{k-1}f')=-\left[-M_{21}\xi_1+\ldots+(-1)^{k-1}M_{k1}\xi_{k-1}\right].\eea
	\par We will show that $f$ has at most finitely many zeros. On the contrary, suppose that $f$ has infinitely many zeros. So we can consider such a point $z_5$ such that
	$f(z_5) = 0$ but $z_5$ is neither a zero nor a pole of $M_{j1}$ $(j=1,2,\ldots,k)$. \par
	Now let, $M_{k1}\not\equiv 0$. Notice that, $\xi_{k-1}$ contain one term, in which, power of $f$ is $0$. Thus we can deduce that $f'(z_5)=0$, which implies that $z_5$ is a multiple zero of $f$ with multiplicity say $p_5\geq2$. Then $z_5$ will be a zero of the left and right hand side of (\ref{e4.44}) with multiplicity $kp_5-1$ and $k(p_5-1)$, respectively, which is a contradiction.\par
	Next, assume $M_{k1}\equiv 0$. If $z_5$ is a simple zero, then $z_5$ is a zero with multiplicity $k-1$ of left hand side of (\ref{e4.44}) and a zero with multiplicity $1$ of right hand side of (\ref{e4.44}), a contradiction. If $z_5$ is a multiple zero with multiplicity $q_5\geq2$, then $z_5$ will respectively be a zero of multiplicity $kq_5-1$ and $k(q_5-1)+1$ of the left and right hand side of (\ref{e4.44}) respectively to yield a contradiction.\par
	Thus, $f$ has at most finitely many zeros.\\	
	\textbf{Sub-case IIC.1.2:} Let $R_5(z)\not\equiv 0$, then (\ref{e4.11}) changes to $f^{n-k-2}(fR_5)=Q_5$ and similarly we get a contradiction.\\	
	\textbf{Sub-case IIC.2:} When $D_0\not\equiv 0$. Proceed similar as {\it Sub-case IB.2} and in this case (\ref{e4.16}), (\ref{e4.17}) and (\ref{e4.18}) changes respectively to \bea\label{e4.45} A_1(f^n+af^{n-2}f')+ \cdots+A_{k+1}(f^n+af^{n-2}f')^{(k)}=Q_2,\eea
	\bea\label{e4.46} \psi_t=\frac{(f^n+af^{n-2}f')^{(t)}}{f^{n-k-2}},\eea $t=1,2,\ldots,k$; $k\geq 3$ and \bea\label{e4.47} f^{n-k-2}R_6=Q_6.\eea Here also $R_6$ is rational. Now, we distinguish following two cases:\\	
	\textbf{Sub-case IIC.2.1:} Let $R_6(z)\equiv 0$, here (\ref{e4.20}) changes to
	\bea\label{e4.48} (f^{k+2}+f^kf')A_1=-(A_2\psi_1+\ldots+A_{k+1}\psi_k).\eea
	We have $f$ has only finitely many zeros. If not, suppose that $f$ has infinitely many zeros. Consider a point $z_6$ such that $f(z_6) = 0$ but $z_6$ is not a zero or a pole of $A_{j}$ $(j=1,2,\ldots,k+1)$. \par
	Now let, $A_{k+1}\not\equiv 0$. Noticing the fact that $\psi_{k}$ contains a term independent of $f$, we can deduce $f'(z_6)=0$, which implies that $z_6$ is a multiple zero of $f$ with multiplicity say $p_6\geq2$. Clearly $z_6$ is a zero of multiplicity $(k+1)p_6-1$ and $(k+1)(p_6-1)$ respectively of the left hand side and right hand side of (\ref{e4.48}), a contradiction.\par
	Next, assume $A_{k+1}\equiv 0$. If $z_6$ is a simple zero, then $z_6$ is a zero with multiplicity $k$ of left hand side of (\ref{e4.48}) and a zero with multiplicity $1$ of right hand side of (\ref{e4.48}), a contradiction. If $z_6$ is a multiple zero with multiplicity $q_6\geq2$, then $z_6$ is a zero of left hand side of (\ref{e4.48}) with multiplicity $(k+1)q_6-1$ and a zero of right hand side of (\ref{e4.48}) with multiplicity $(k+1)q_6-k$, again gives a contradiction. It follows that $f$ has at most finitely many zeros.\\	
	\textbf{Sub-case IIC.2.2:} If $R_6(z)\not\equiv 0$, then (\ref{e4.21}) changes to $$f^{n-k-3}(fR_6)=Q_6$$ and adopting similar procedure we can get a contradiction. \\
	
	As usual from the above two cases we conclude that $f$ is a transcendental meromorphic function with finite zeros and poles. Now by {\it{Lemma \ref{l3.4}}}, we can say that
	\bea\label{e4.49}f(z) = q(z)e^{P(z)},\eea where $q(z)$ is a non-vanishing rational function, and $P(z)$ is a non-constant polynomial. Substituting (\ref{e4.49}) into (\ref{e2.1}) yields
	\bea\label{e4.50} q^n(z)e^{nP(z)}+aq^{n-2}(q'+qP')e^{(n-1)P(z)}+\sum\limits_{j=0}^{d}\varphi_j(z)e^{jP(z)}=\sum\limits_{i=1}^{k}p_i(z)e^{\alpha_i(z)},\eea
	where $\varphi_j(z)$ are rational functions.\\
	
	Now we can rearrange $\{1,2,\ldots,d\}$ to $\{\sigma_1,\sigma_2,\ldots,\sigma_d\}$ such that $\varphi_j(z)e^{jP(z)}=\varphi_{\sigma_j}(z)e^{{\sigma_j}P(z)}$ for $j=1,\ldots,d$ and $\{1,2,\ldots,k\}$ to $\{\mu,\nu,\kappa_1,\kappa_2\ldots,\kappa_{k-2}\}$ such that $p_1(z)e^{\alpha_1(z)}=p_{\mu}(z)e^{\alpha_{\mu}(z)}$, $p_2(z)e^{\alpha_2(z)}=p_{\nu}(z)e^{\alpha_{\nu}(z)}$ and $p_i(z)e^{\alpha_i(z)}=p_{\kappa_{i-2}}(z)e^{\alpha_{\kappa_{i-2}}(z)}$ for $i=3,4,\ldots,k$. Then (\ref{e4.50}) can be written as \bea\label{e4.51}&& q^n(z)e^{nP(z)}+aq^{n-2}(q'+qP')e^{(n-1)P(z)}+\sum\limits_{j=0}^{d}\beta_{\sigma_j}(z)e^{{\sigma_j}P(z)}\nonumber\\&&=p_{\mu}(z)e^{\alpha_{\mu}(z)}+p_{\nu}(z)e^{\alpha_{\nu}(z)}+\sum\limits_{i=3}^{k}p_{\kappa_{i-2}}(z)e^{\alpha_{\kappa_{i-2}}(z)}.\eea 
	
	Since $n>n-1>d$, $\deg(\alpha_{\kappa_i}-\alpha_{\kappa_j})\geq 1$ $(1\leq \kappa_i\neq \kappa_j\leq k)$ and $q$, $p_{\kappa_{i-2}}$ $(i=3,4,\ldots,k)$ are all non-zero rational functions and $q'+qP'\not\equiv 0$, using {\it {Lemma \ref{l3.3}}} on (\ref{e4.51}) we have
	\beas nP=\alpha_{\mu}+\mathcal{C}_{\mu},(n-1)P=\alpha_{\nu}+\mathcal{C}_{\nu},\sigma_1P=\alpha_{\kappa_1}+\mathcal{C}_1,\sigma_2P=\alpha_{\kappa_2}+\mathcal{C}_2,\ldots,\sigma_{k-2}P=\alpha_{\kappa_{k-2}}+\mathcal{C}_{k-2},\eeas where $\mathcal{C}_i$ $(i=1,2,\ldots,k-2)$ are constants such that \beas q^n=e^{-\mathcal{C}_\mu}p_{\mu},aq^{n-2}(q'+qP')=e^{-\mathcal{C}_\nu}p_{\nu},\beta_{\sigma_1}=e^{-\mathcal{C}_1}p_{\kappa_1},\beta_{\sigma_2}=e^{-\mathcal{C}_2}p_{\kappa_2},\ldots,\beta_{\sigma_{k-2}}=e^{-\mathcal{C}_{k-2}}p_{\kappa_{k-2}}.\eeas
	Also, $\beta_0(z)=0$ and $\beta_{\sigma_j}(z)\equiv 0$ for all $\sigma_j\neq \sigma_1,\sigma_2,\ldots,\sigma_{k-2}$ with $0\leq \sigma_j\leq d$.\par Therefore, \beas \alpha'_{\mu}:\alpha'_{\nu}:\alpha'_{\kappa_1}:\ldots:\alpha'_{\kappa_{k-2}}=n:n-1:\sigma_1:\ldots:\sigma_{k-2},\eeas
	\beas nP'=\alpha'_{\mu},\;(n-1)P'=\alpha'_{\nu}\;and\;P_d(z,f)=\sum\limits_{i=1}^{k-2}p_{\kappa_i}(z)e^{\alpha_{\kappa_i}(z)}.\eeas

	This completes the proof of the theorem.
\end{proof}

\end{document}